\documentclass[11pt]{article}
\usepackage[T1]{fontenc}
\usepackage{lmodern}
\usepackage[a4paper,margin=26mm]{geometry}
\usepackage{amsmath,amssymb,amsthm,mathtools}
\usepackage{booktabs,tabularx,tikz}
\usepackage{microtype,xurl}
\usepackage{caption}
\usepackage{titlesec}
\titleformat{\section}{\normalfont\Large\bfseries}{\thesection}{1em}{}
\titleformat{\subsection}{\normalfont\large\bfseries}{\thesubsection}{1em}{}
\titlespacing*{\section}{0pt}{2.4ex plus .5ex minus .2ex}{1.2ex plus .2ex}
\titlespacing*{\subsection}{0pt}{2ex plus .4ex minus .2ex}{1ex plus .2ex}
\titlespacing*{\paragraph}{0pt}{2ex plus .3ex minus .2ex}{1em}
\usepackage[numbers,sort&compress]{natbib}
\usepackage[hidelinks]{hyperref}
\providecommand{\doi}[1]{\href{https://doi.org/#1}{\nolinkurl{doi:#1}}}
\newtheorem{theorem}{Theorem}[section]
\newtheorem{lemma}[theorem]{Lemma}
\newtheorem{proposition}[theorem]{Proposition}
\theoremstyle{definition}
\newtheorem{example}[theorem]{Example}
\theoremstyle{plain}

\newcommand{\D}{\mathbb D}
\newcommand{\HH}{\mathbb H}
\newcommand{\C}{\mathbb C}
\newcommand{\Om}{\Omega}
\newcommand{\re}{\operatorname{Re}}
\newcommand{\im}{\operatorname{Im}}
\newcommand{\tr}{\operatorname{tr}}
\newcommand{\arcosh}{\operatorname{arcosh}}
\newcommand{\atanh}{\operatorname{artanh}}
\newcommand{\dd}{\mathrm{d}}

\title{The Exact Stabilizability Threshold\\for the Belgian Chocolate Problem}
\author{
  Xiang Ji\textsuperscript{1}, Long Wang\textsuperscript{2}, Wensheng Yu\textsuperscript{1}\\[0.55em]
  {\normalsize \textsuperscript{1}\,School of Electronic Engineering,}\\
  {\normalsize Beijing University of Posts and Telecommunications, Beijing 100876, China}\\[0.15em]
  {\normalsize \textsuperscript{2}\,Center for Systems and Control, College of Engineering,}\\
  {\normalsize Peking University, Beijing 100871, China}
}
\date{}
\hypersetup{pdfauthor={Xiang Ji, Long Wang, Wensheng Yu},pdftitle={The Exact Stabilizability Threshold for the Belgian Chocolate Problem}}

\begin{document}
\maketitle
\vspace{-1.4em}

\begin{abstract}
The Belgian chocolate problem arose from Blondel's work on simultaneous stabilization of linear systems. Its generalized form asks which parameters of a second-order plant admit a proper, stable, minimum-phase rational controller. Early work raised the known lower bound on the stabilizability threshold \(\delta_*\) by constructing stabilizing controllers and lowered its upper bound using complex analysis. We construct a balanced slit domain and compare harmonic primitives to prove that the feasible positive parameters form the open interval \((0,\delta_*)\), where \(\delta_*=\cos(\pi d_*/2)\) and \(d_*\) is the harmonic measure at the origin of the lower semicircle in the balanced domain. At balance the curved arc has the same harmonic measure as the two slit banks combined. No finite-order controller stabilizes the plant at \(\delta_*\). The corresponding disk interpolation problem has exactly two solutions at its critical parameter. For controllers with sixth-degree denominators, allowing a sixth-degree numerator enlarges the feasible parameter set. With a constant numerator, the feasible parameters can approach the continued-fraction bound arbitrarily closely at denominator degrees one through four. At every fixed degree five or higher, the feasible parameters stay a positive distance below that bound. We verify the threshold theorems and both fixed-degree results in Lean~4.
\end{abstract}


\section{Introduction and main result}\label{sec:intro}

Simultaneous stabilization asks whether a single controller can stabilize every plant in a given family. For a single plant, strong stabilization requires the controller itself to be stable \cite{Youla}. The reduction of a three-plant problem leads to a one-plant problem with the additional requirement that the controller be minimum-phase \cite[pp.~149--150]{Blondel}. Blondel's Belgian chocolate problem asks whether a proper, stable, minimum-phase rational controller exists at \(\delta=0.9\) for the second-order plant
\[
 P_\delta(s)=\frac{b(s)}{a_\delta(s)},\qquad
 a_\delta(s)=s^2-2\delta s+1,\qquad b(s)=s^2-1.
\]
The generalized problem asks for the full set of parameters \(\delta\) for which such a controller exists. We use the polynomial formulation given in \cite[Section~5]{GuanHeWangYu} and \cite[Eq.~(4)]{LiWangYu}.

A real polynomial is \emph{strictly Hurwitz stable} if it is nonzero and all its complex roots have strictly negative real part. Write a controller as \(c(s)=y(s)/x(s)\). Strict Hurwitz stability of \(x\) and \(y\) places the poles and finite zeros of \(c\), respectively, in the open left half-plane, giving stability and minimum phase. The controller is proper when \(\deg y\leq\deg x\). A parameter \(\delta\) is \emph{feasible} if such \(x,y\) make the following closed-loop polynomial strictly Hurwitz stable
\begin{equation}
 Z(s)=a_\delta(s)x(s)+b(s)y(s)
 \label{eq:poly}
\end{equation}
Write \(\delta_*\) for the supremum of the feasible positive parameters.

At \(\delta=0\), the choice \(x(s)=s+1\), \(y(s)=1/2\) gives \(Z(s)=(s+1)(s^2+s/2+1/2)\), which is strictly Hurwitz stable. At \(\delta=1\), both \(a_\delta\) and \(b\) vanish at \(s=1\), so \(Z(1)=0\). If \(\delta>1\), the two positive roots \(r_-<1<r_+\) of \(a_\delta\) give \(Z(r_\pm)=b(r_\pm)y(r_\pm)\) with opposite signs, since \(y\) has no positive root. Thus \(Z\) has a positive root. The positive-parameter analysis therefore reduces to \(0<\delta<1\).

The existence of a threshold separating feasible from infeasible positive parameters follows from the early work of Blondel and collaborators. Estimates by Rupp and by Blondel and collaborators gave the lower bound \(\delta_*>0.761594\) \cite[Section~5]{GuanHeWangYu}. In 2002, Patel, Deodhare, and Viswanath answered Blondel's original question by finding an eleventh-order controller for \(\delta=0.9\) by randomized search. They also improved the lower bound to \(0.937\) \cite{Patel}. Further controller constructions raised it to \(0.96\) in 2006 by nonsmooth optimization and to \(0.973974\) in 2007 by global optimization \cite{BHLO,ChangSahinidis}. In the same year, Batra's estimates for the locations of zeros and one-points of holomorphic functions yielded an upper bound of about \(0.99999712\) \cite{Batra2007,BE}. Polynomial inequalities also gave bounds for controllers of prescribed order in 2008 \cite{HeWY}.

Eremenko discussed the connection with holomorphic avoidance in 2012 \cite{Eremenko}. This connection led to an extremal interpolation problem on the disk. In 2013, Bergweiler and Eremenko used hyperbolic geodesic lengths in that formulation to prove \(\delta_*<0.999579\) \cite[Theorem~11.1]{BE}. In 2014, Li and Yu refined the Hempel--Smith geodesic-length estimate by using elliptic modular functions and obtained \(r\approx0.01328891154\) for the reciprocal radius in their calculation \cite[Eq.~(10)]{LiYu}. In 2018, the algebraic structure of limiting controllers allowed Charles and Boston to establish feasibility at \(\delta=0.9808348\), narrowing the known range to \(0.9808348<\delta_*<0.999579\) \cite[Theorem~5.1]{CB}.

In September 2026, Cheng extended the certified feasible range by giving an explicit rational polynomial certificate at \(\delta=0.9846316110433085\) \cite{ChengBCP}. In November 2026, Liu combined harmonic measure with a Cauchy integral argument to characterize the exact threshold and obtained \(\delta_*\approx0.987237027\) by solving the balance equation numerically \cite[Theorem~1.1 and Sections~5--6]{Liu}.

We prove the exact stabilizability threshold by comparing harmonic primitives (Section~\ref{sec:comparison}) and classify all disk interpolants at the critical parameter. We also resolve two conjectures for controllers of prescribed order \cite{HeWY}. Section~\ref{sec:fixed-degree} gives a sixth-order counterexample and determines when the constant-numerator bound is sharp. The formula for the global threshold is determined by a balance condition on the boundary of the following slit domain.

Write \(\C\) for the complex plane, \(i\) for the imaginary unit, and \(\D=\{z\in\C: |z|<1\}\) for the open unit disk. Consider the curved arc and the lower unit semicircle
\[
 C=\{\exp(-\sqrt{\theta(\pi-\theta)}+i\theta):0<\theta<\pi\},
 \qquad L=\{e^{i\theta}:\pi<\theta<2\pi\}.
\]
The arcs join at \(-1\) and \(1\) to form a simple closed curve. Let \(G\) be the region it encloses, a Jordan domain. For \(-1<\alpha<0\), remove the real segment \([-1,\alpha]\) from \(G\). This gives the slit domain \(D_\alpha=G\setminus[-1,\alpha]\) with tip \(\alpha\).

At an interior point of the slit, approach from above or below corresponds to two distinct points of the prime-end boundary, the boundary used for conformal correspondence \cite[Section~2.4]{Pommerenke}. Let \(S_\alpha\) denote the union of these two banks in that boundary. Take a conformal map from \(\D\) onto \(D_\alpha\) that sends \(0\) to \(0\). The harmonic measure \(\omega_{D_\alpha}(0,E)\) of a boundary portion \(E\) is the fraction of the unit circle that corresponds to \(E\) under the prime-end boundary correspondence. This fraction is independent of the choice of the centered map, since any two such maps differ by a rotation of \(\D\) \cite[Section~4.4]{Pommerenke}. In particular, \(S_\alpha\) contributes the sum of the weights of the two banks. The slit is \emph{balanced} when this sum equals the weight of \(C\).

To compare the harmonic measures as the slit tip varies, write
\[
 a(\alpha)=\omega_{D_\alpha}(0,C),\qquad
 b_0(\alpha)=\omega_{D_\alpha}(0,S_\alpha),\qquad
 d(\alpha)=\omega_{D_\alpha}(0,L).
\]
The balance condition is \(a(\alpha)=b_0(\alpha)\).

\begin{theorem}[Stabilizability threshold]\label{th:main}
There is a unique balanced slit tip \(\alpha_*\in(-1,0)\), characterized by
\[
 \omega_{D_{\alpha_*}}(0,C)=\omega_{D_{\alpha_*}}(0,S_{\alpha_*}).
\]
Let \(d_*=d(\alpha_*)\) be the harmonic measure of the lower semicircle at balance. The feasible positive parameters for \eqref{eq:poly} form the interval \((0,\delta_*)\), where
\begin{equation}
 \delta_*=\cos\!\left(\frac{\pi d_*}{2}\right).
 \label{eq:delta}
\end{equation}
For the parameter domain \([0,1]\), the feasible set is \([0,\delta_*)\). The threshold is an unattained supremum over finite controller orders.
\end{theorem}

Our Lean~4 formalization starts from the concrete slit domains and the feasibility definitions and establishes the exact threshold formulas, both feasible intervals, and their endpoint conclusions. The main theorem of Liu's public Lean development assumes certain properties of the harmonic-measure data, the existence of a critical interpolant, and the existence of a half-disk comparison function \cite{LiuLean}. Applying that theorem to the original problem in Lean requires proving that these hypotheses hold for the concrete slit domains.

\section{The disk interpolation problem}\label{sec:disk}

For \(0<t<1\), a \emph{disk interpolant at \(t\)} is a holomorphic function \(f\) on the unit disk with the real symmetry
\begin{equation}\label{eq:disk}
 f(\bar z)=\overline{f(z)}\qquad(z\in\D).
\end{equation}
It takes the value \(0\) only at the origin and the value \(1\) only at \(it\) and \(-it\), with nonzero derivative at all three points. Here \(\bar z\) denotes complex conjugation of \(z\). The symmetry makes \(f'(0)\) real. Replacing \(f(z)\) by \(f(-z)\), when necessary, gives the normalization \(f'(0)>0\).

\begin{theorem}[Disk threshold and uniqueness]\label{th:disk}
Let \(\alpha_*,d_*\) be as in Theorem~\ref{th:main}, and put
\[
 t_*=\tan\!\left(\frac{\pi d_*}{4}\right).
\]
The set of \(t\in(0,1)\) admitting a disk interpolant is exactly \([t_*,1)\). At \(t_*\), there are exactly two interpolants, \(f_*\) and \(z\mapsto f_*(-z)\), with \(f_*'(0)>0\).
\end{theorem}

Let \(\delta\) be feasible, with strictly Hurwitz polynomials \(x,y,Z\) satisfying \eqref{eq:poly}. The reduction to the disk uses the Cayley transform \cite[Sections~10--11]{BE}.
Set \(s=(1+z)/(1-z)\) and
\begin{equation}
 R(z)=\frac{(s^2-1)y(s)}{Z(s)},\qquad
 1-R(z)=\frac{(s^2-2\delta s+1)x(s)}{Z(s)} .
 \label{eq:cayley}
\end{equation}
The Cayley map sends \(\D\) onto the open right half-plane, where \(x,y,Z\) have no zeros. Thus \(R\) is holomorphic in \(\D\). Its only zero corresponds to \(s=1\), and its two \(1\)-points correspond to \(s=\delta\pm i\sqrt{1-\delta^2}\). All are simple, and the real coefficients give the symmetry. Hence \(R\) is a rational disk interpolant with
\begin{equation}
 t=\sqrt{\frac{1-\delta}{1+\delta}},\qquad
 \delta=\frac{1-t^2}{1+t^2}.
 \label{eq:parameter}
\end{equation}

The proof first constructs a candidate from the balanced slit in Section~\ref{sec:balance}. In Section~\ref{sec:comparison}, the annuli fix the boundary markings of a competing lift, and the harmonic comparison gives the threshold and endpoint classification. Section~\ref{sec:finite} translates these conclusions back to finite-order controllers.

\section{The balanced slit and the candidate}\label{sec:balance}

\subsection{The balance equation}

Let \(\Om=\{z\in\D:\im z>0\}\) be the upper half-disk. Use \(z\) for its coordinate and \(w\) for the coordinate on \(D_\alpha\). The point \(-1\) has two prime ends. Denote by \(P\) the one between \(C\) and the upper bank, and by \(Q\) the one between the lower bank and \(L\). The branch \(w\mapsto\sqrt{w-\alpha}\) with positive real part sends \(D_\alpha\) to a Jordan domain. Its boundary consists of the two unfolded exterior arcs and the segment between \(-i\sqrt{1+\alpha}\) and \(i\sqrt{1+\alpha}\). The exterior arcs meet at the image of \(1\). Their other endpoints are the distinct images of \(P,Q\), and the slit tip maps to the midpoint of the segment. The boundary correspondence theorem for Jordan domains, or the prime-end boundary correspondence in \(D_\alpha\), gives a conformal map \(\Phi_\alpha:\Om\to D_\alpha\) marked by
\[
 -1,\ 0,\ 1\ \longmapsto\ 1,\ P,\ Q
\]
in positive boundary order \cite[Chapter~2]{Pommerenke}. Thus \((-1,0)\) corresponds to \(C\), \((0,1)\) to the two slit banks, and the upper semicircle to \(L\). The three weights are positive and sum to one.

\begin{figure}[htbp]
\centering
\begin{tikzpicture}[x=1.85cm,y=1.55cm,font=\small]
\begin{scope}
\draw[thick,blue!65!black] (-1,0) arc[start angle=180,end angle=0,radius=1];
\draw[thick] (-1,0)--(1,0);
\fill (0,.39) circle (1.4pt) node[right] {\(it_*\)};
\node[below] at (-1,0){\(-1\)};
\node[below] at (0,0){\(0\)};
\node[below] at (1,0){\(1\)};
\node at (0,.73){\(\Omega\)};
\end{scope}
\draw[->,thick] (1.29,.23)--(2.04,.23) node[midway,above] {\(\Phi_*\)};
\begin{scope}[shift={(3.45,.45)}]
\draw[thick,blue!65!black] (-1,0) arc[start angle=180,end angle=360,radius=1];
\draw[thick,samples=151,domain=0:180,variable=\ang]
  plot ({exp(-sqrt(max(0,(\ang*pi/180)*(pi-\ang*pi/180))))*cos(\ang)},
        {exp(-sqrt(max(0,(\ang*pi/180)*(pi-\ang*pi/180))))*sin(\ang)});
\draw[very thick] (-1,0)--(-.34,0) node[below] {\(\alpha_*\)};
\fill (0,0) circle (1.4pt) node[below right] {\(0\)};
\node[above left] at (-1,0){\(P,Q\)};
\node[right] at (1,0){\(1\)};
\node[above] at (0,.22){\(C\)};
\node at (.20,-.82){\(L\)};
\end{scope}
\end{tikzpicture}
\caption{The two halves of the diameter and the upper semicircle correspond to \(C\), the two slit banks, and \(L\), respectively. The marked interior point and slit tip are shown schematically.}
\label{fig:slit}
\end{figure}
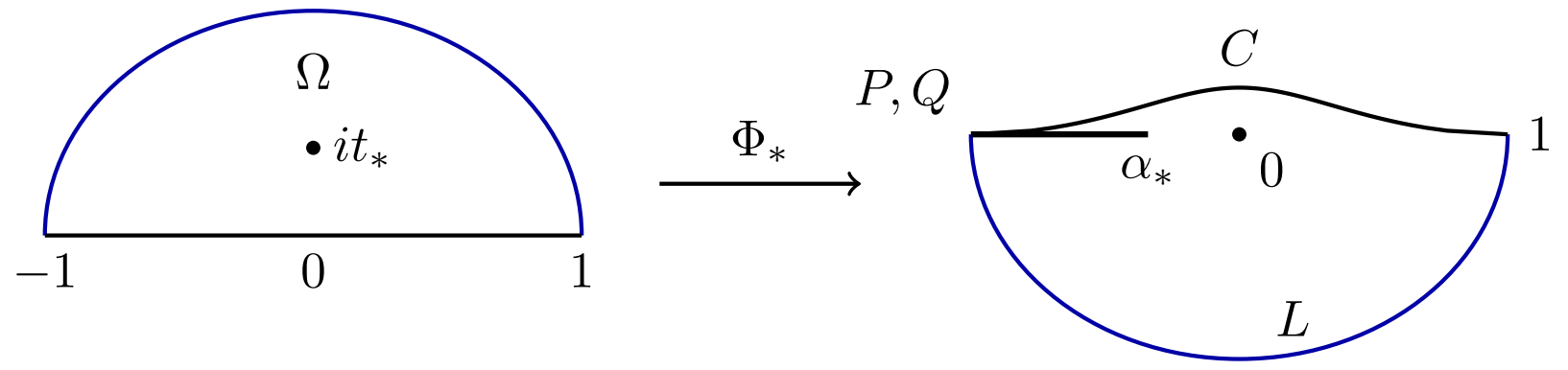

\begin{proposition}\label{pr:balance}
There is exactly one slit tip \(\alpha\) in \((-1,0)\) for which the curved arc and the two slit banks have equal harmonic measures at the origin, that is, \(a(\alpha)=b_0(\alpha)\).
\end{proposition}

\begin{proof}
For \(\alpha<\beta\), restrict the harmonic measure functions for \(D_\alpha\) to \(D_\beta\). On the added slit segment \((\alpha,\beta)\), the old harmonic measure of \(C\) is strictly positive and that of the old slit is strictly less than one. On the common boundary the corresponding values agree except at finitely many prime ends. Boundary comparison gives \(a(\beta)<a(\alpha)\) and \(b_0(\beta)>b_0(\alpha)\), since the added segment has positive harmonic measure in \(D_\beta\).

If a sequence of slit tips \((\alpha_n)\) converges to a point \(\alpha\) in \((-1,0)\), the pointed domains \((D_{\alpha_n},0)\) converge to \((D_\alpha,0)\) in the Carath\'eodory kernel sense. For these simply connected domains, this means locally uniform convergence of the Riemann maps from \(\D\), normalized to send \(0\) to \(0\) with positive derivative there. The harmonic measures on the planar boundaries, obtained by assigning each prime end to its corresponding point on the planar boundary, converge weakly. That is, their integrals against continuous test functions converge \cite[Section~2.1 and Theorem~2.3]{BRY}. In the common compact set \(\overline C\cup\overline L\cup[-1,0]\), the relative boundary of either \(C\) or \(L\) consists of the points \(\pm1\). These represent only finitely many prime ends and have harmonic measure zero. Thus \(a,d\) are continuous. The identity \(b_0=1-a-d\) gives continuity of \(b_0\).

As \(\alpha\) approaches \(-1\) from the right, put \(\epsilon=1+\alpha\). The slit lies in \(|w+1|\le\epsilon\). Use the harmonic function \(\log(2/|w+1|)/\log(2/\epsilon)\) as a barrier. It is at least \(1\) on the slit and nonnegative on the outer boundary, so it bounds the harmonic measure of the slit from above
\[
 b_0(\alpha)\le\frac{\log 2}{\log(2/\epsilon)}.
\]
This upper bound tends to zero as \(\epsilon\) tends to zero. The loss of \(C\)-measure relative to \(G\) lies between \(0\) and \(b_0(\alpha)\), by comparing the two harmonic measure functions on the slit and the common exterior boundary. Hence \(a(\alpha)\) tends to the positive \(C\)-measure in \(G\).

As \(\alpha\) approaches \(0\) from the left, put \(\epsilon=-\alpha\) and fix \(0<r<e^{-\pi/2}\), with \(\epsilon<r/2\). The disk of radius \(r\) lies in \(G\). In the Brownian exit interpretation of harmonic measure \cite[Section~1.1]{BRY}, a path from \(0\) to either exterior arc before hitting the slit must first leave the disk of radius \(r\) without crossing the slit. In the coordinate \(\sqrt{w+\epsilon}\), the image starts at \(\sqrt\epsilon\) in the right half-plane and must first reach the right semicircle of radius \(\sqrt{r-\epsilon}\). For the half-plane computation, use the Joukowski map \(J(z)=z+1/z\), which sends \(\Om\) conformally onto the lower half-plane. Scale the right half-disk to unit radius and rotate it by multiplication by \(i\) before applying \(J\). The half-plane Poisson kernel at \(x+iy\), as a density in the real boundary coordinate \(u\), is \(|y|/(\pi((u-x)^2+y^2))\) \cite[Chapter~7]{HFT}. Integrating it over the image of the semicircle gives
\[
 1-b_0(\alpha)\le
 \frac4\pi\arctan\sqrt{\frac{\epsilon}{r-\epsilon}}.
\]
As \(\epsilon\) tends to zero, the upper bound tends to zero, so \(b_0\) tends to \(1\) and \(a\) tends to zero. The limits at the two ends give opposite signs for \(a-b_0\). Continuity and strict decrease then give a unique solution.
\end{proof}

Set \(\Phi_*=\Phi_{\alpha_*}\) and \(p=\Phi_*^{-1}(0)\). Under \(J\), the two pieces of the real diameter map to the rays \((-\infty,-2)\) and \((2,\infty)\). Their harmonic measures at \(J(p)\) coincide exactly when \(\re J(p)=0\), by the half-plane Poisson kernel. Since \(\re J(p)=(1+|p|^{-2})\re p\), the balance condition gives \(p=it_0\), where \(t_0\) is the positive imaginary part of \(p\). The upper semicircle maps to \([-2,2]\). Integrating the same kernel at \(J(it_0)=-i(1-t_0^2)/t_0\) gives
\begin{equation}\label{eq:mass}
 d_*=\frac2\pi\arctan\!\left(\frac{2t_0}{1-t_0^2}\right)
     =\frac4\pi\arctan t_0 .
\end{equation}
It follows that \(t_0=\tan(\pi d_*/4)=t_*\), and \(\cos(\pi d_*/2)=(1-t_*^2)/(1+t_*^2)\).

A half-turn maps Liu's positive slit to the slit used here, sending his slit tip to \(\alpha\) and matching his curved and outer arcs with \(C\) and \(L\). The three weights sum to one, so Liu's balance condition is equivalent to \(a(\alpha)=b_0(\alpha)\) \cite[Theorem~1.1 and Section~3]{Liu}.

\subsection{The modular covering and the candidate}

Let \(X=\C\setminus\{0,1\}\). The modular function \(\lambda:\HH\to X\) is the universal cover with deck group \(\Gamma(2)/\{\pm I\}\), identifying matrices \(M\) and \(-M\) that act by the same fractional linear transformation. Here \(\HH=\{\tau\in\C:\im\tau>0\}\) is the upper half-plane and \(\Gamma(2)=\{M\in\mathrm{SL}_2(\mathbb Z):M\equiv I\pmod 2\}\) acts on it by fractional linear transformations. \(\mathrm{SL}_2(\mathbb Z)\) is the group of integral \(2\times2\) matrices of determinant one and \(I\) is its identity. The source of this covering is simply connected, hence the name universal cover. Its deck transformations are the biholomorphisms \(M\) of \(\HH\) with \(\lambda\circ M=\lambda\) \cite[Section~2]{BE}. At the cusp \(\infty\), the translation \(\tau\mapsto\tau+2\) preserves \(\lambda\). Put \(q=e^{\pi i\tau}\). Since \(\lambda\) has period \(2\), its value depends only on \(q\). We write \(\lambda(\tau)=\Lambda(q)\) for the same modular function expressed in this coordinate. As the imaginary part of \(\tau\) tends to infinity, \(q\) tends to \(0\). The function \(\Lambda\) extends holomorphically to \(q=0\). Its Taylor series starts with \(16q\), and the remainder \(O(q^2)\) starts at degree two. We use the standard relations \cite[Section~2]{BE} and \cite[Sections~23.15, 23.17, and 23.18]{DLMF}
\begin{gather}
 \Lambda(q)=16q+O(q^2),\qquad
 \lambda(-1/\tau)=1-\lambda(\tau),\qquad
 \lambda(\tau+1)=\frac{\lambda(\tau)}{\lambda(\tau)-1},
 \label{eq:lambda}\\
 \lambda(-\bar\tau)=\overline{\lambda(\tau)},\qquad
 0<\lambda(iy)<1\quad(y>0).
 \nonumber
\end{gather}
The coordinate \(w=-\exp(-\pi i/\tau)\) identifies the quotient of \(\HH\) by the translation \(-1/\tau\mapsto-1/\tau+2\) with \(\D\setminus\{0\}\). The induced map
\[
 U(w)=1-\Lambda(-w)
\]
is a covering onto \(X\) on the punctured disk and satisfies \(U(w)=\lambda(\tau)\). At the puncture it extends with \(U(0)=1\) and \(U'(0)=16\). Thus \(U\) has no disk zeros and takes the value \(1\) only at \(0\), where its derivative is nonzero.

For \(\tau=iy\), the coordinate \(w\) runs through \((-1,0)\) and \(0<U(w)<1\). For \(\tau=1+iy\),
\[
 \theta=\arg w=\frac{\pi y^2}{1+y^2},\qquad
 -\log|w|=\frac{\pi y}{1+y^2}=\sqrt{\theta(\pi-\theta)}.
\]
It follows that this vertical line maps onto \(C\), where \(U\) takes negative real values by \eqref{eq:lambda}.

\begin{proposition}\label{pr:candidate}
The function \(U\circ\Phi_*\) on \(\Om\) extends to a disk interpolant \(f_*\) at \(t_*\) with \(f_*'(0)>0\).
\end{proposition}

\begin{proof}
In \(\Om\), the composition has no zeros and has the single simple \(1\)-point \(it_*\). On each open real segment it has continuous real boundary values, negative on \((-1,0)\) and in \((0,1)\) on \((0,1)\). The two values from the banks coalesce at the slit tip, where \(U\) is continuous. Schwarz reflection therefore applies along both open segments.

To determine the order of vanishing at the source point \(0\), examine the image near \(P\). Choose the branch of \(\log(-w)\) tending to zero as \(w\) approaches \(-1\) and set \(T=\log(-w)/(\pi i)\) and let \(\tau=-1/T\). The bordering arcs of \(D_{\alpha_*}\) are the upper slit bank and \(C\). They map to the lines \(\re\tau=0,1\). The interior adjacent to \(P\) maps between these lines at large positive imaginary part. In fact, the branch is in \(\HH\), and near the upper bank a displacement into the domain gives \(\re\tau>0\). The bordering curved arc is precisely the image of \(\re\tau=1\) computed above. As \(w\) approaches \(P\) inside the domain, the real part of \(\tau\) lies between \(0\) and \(1\) and its imaginary part tends to infinity. Thus \(q=e^{\pi i\tau}\) approaches \(0\) from the upper half-plane. The reflected function extends across \(0\) by removability with value zero. The real Taylor coefficients of \(\Lambda\) give
\[
 \frac{\im\Lambda(q)}{\im q}=16+O(|q|)\qquad(\im q>0),
\]
uniformly as \(q\) approaches \(0\) from the upper half-plane, using \(|\im q^n|\le n|q|^{n-1}\im q\). Thus \(\im f_*>0\) in a small upper half-neighborhood of \(0\). A real holomorphic function with first nonzero term \(Az^m\) and \(A>0\) has a negative imaginary part along some ray in that half-neighborhood if \(m\ge2\). Here \(A>0\) because the positive real values are positive. Hence \(m=1\). Reflection gives the other simple \(1\)-point at \(-it_*\), and the real-segment values exclude any other \(0\)- or \(1\)-points.
\end{proof}

\section{Global comparison and the extremal interpolant}\label{sec:comparison}

\subsection{Comparison of harmonic primitives}

For a bounded harmonic function \(v\) on \(\D\), write \(v^*(\xi)\) for the radial boundary value at a point \(\xi\) of the unit circle. It is the limit of \(v(r\xi)\) along the radius from inside the disk, and exists almost everywhere with respect to arc length \cite[Chapter~6]{HFT}. On arcs of continuous extension we use the same superscript for the boundary restriction. A nonnegative essential lower limit along an arc at an endpoint means that for every \(\varepsilon>0\) the function is at least \(-\varepsilon\) almost everywhere on a sufficiently short portion of that arc next to the endpoint.

\begin{lemma}[Boundary comparison]\label{lem:primitive}
Let \(W_j=u_j+iv_j\) be holomorphic in \(\D\), with \(W_j(0)=0\) and \(v_j\) bounded for \(j=0,1\). Partition the unit circle into three open arcs and their three endpoints. On the first arc both functions extend continuously and satisfy \(v_j=g(u_j)\) for one single-valued function \(g\). On the other two arcs, the radial boundary values satisfy \(v_1^*-v_0^*\ge0\) almost everywhere. At each endpoint the essential lower limits of this difference along the two adjacent arcs are nonnegative. Then \(W_1=W_0\).
\end{lemma}

\begin{proof}
For \(W=u+iv\) as above define the harmonic primitive
\[
 H(z)=\im\!\int_0^z \frac{W(\zeta)}{\zeta}\,\dd\zeta .
\]
The integrand has a removable singularity at \(0\). Applying the same formula to \(W_0,W_1\) defines \(H_0,H_1\). In polar coordinates \(z=re^{i\theta}\), the partial derivatives with respect to radius and angle, \(H_r\) and \(H_\theta\), satisfy
\begin{equation}
 rH_r=v,\qquad H_\theta=u .
 \label{eq:derivs}
\end{equation}
The first identity gives
\(\sup_{z\in\D}|H(z)-H(rz)|\le\|v\|_\infty\log(1/r)\).
The dilated functions \(H(rz)\) are continuous on the closed disk and uniformly Cauchy as \(r\) approaches \(1\) from below. Their uniform limit extends \(H\) continuously to the boundary. Integrating the second identity on compact subarcs of the first arc and passing to the boundary shows that the derivative along the circle with respect to the angle is \(u^*\), the continuous boundary value of \(u\).

Let \(D=H_1-H_0\) and \(h=v_1-v_0\). If \(D\) is nonconstant, it has a boundary minimum \(m=D(\xi)\) with \(D(0)>m\). For \(\xi,\zeta\in\partial\D\) and \(0<r<1\), the Poisson kernel with respect to normalized arc length is \((1-r^2)/|\zeta-r\xi|^2\). It is bounded below by \((1-r)/(1+r)\) \cite[Chapter~1]{HFT}. Applied to the nonnegative boundary values of \(D-m\), this gives
\begin{equation}
 \frac{m-D(r\xi)}{1-r}\le -\frac{D(0)-m}{1+r}.
 \label{eq:comparison-upper}
\end{equation}
If \(\xi\) lies on the first arc, vanishing of the tangential derivative gives \(u_1^*(\xi)=u_0^*(\xi)\). The graph condition then gives \(h^*(\xi)=0\), and continuity gives a radial limit of zero for \(h(r\xi)\). Otherwise, for each \(\epsilon>0\) a neighborhood arc of \(\xi\) has \(h^*\ge-\epsilon\) almost everywhere. The Poisson representation for bounded harmonic functions \cite[Chapters~1 and~6]{HFT} and the vanishing Poisson mass of the complementary arc give a nonnegative lower limit of \(h(r\xi)\) as \(r\) approaches \(1\) from below. In both cases, radial integration of \eqref{eq:derivs} gives
\[
 \frac{m-D(r\xi)}{1-r}
   =\frac{1}{1-r}\int_r^1h(\rho\xi)\,\frac{\dd\rho}{\rho}.
\]
The lower limit of the right-hand side is nonnegative because \(\log(1/r)/(1-r)\) tends to \(1\) as \(r\) approaches \(1\) from below. This contradicts \eqref{eq:comparison-upper}. Hence \(D\) is constant. Equation \eqref{eq:derivs} gives \(v_1=v_0\), and holomorphicity and normalization imply \(W_1=W_0\).
\end{proof}

\subsection{The marked lift of a competitor}\label{sec:lift}

A lift of an interpolant \(f\) through \(U\) is a map \(\phi\) such that \(U\circ\phi=f\) on \(\Om\). The boundary markings specify where this lift sends the two open halves of the real diameter.

\begin{proposition}\label{pr:lift}
Let \(0<t\le1/10\) and let \(f\) be a disk interpolant at \(t\) with \(f'(0)>0\). There is a holomorphic map \(\phi:\Om\to\D\) with \(U\circ\phi=f\) and \(\phi(it)=0\). The same lift extends analytically across each open real segment and satisfies
\begin{equation}
 \phi((-1,0))\subset C,\qquad
 \phi((0,1))\subset(-1,0).
 \label{eq:marks}
\end{equation}
\end{proposition}

\begin{proof}
Lift \(f\) locally through \(\lambda\) on \(\D\setminus\{0,it,-it\}\), based at \(z_b=-1/2\). Let \(h=t/5\). For each \(p\in\{0,it,-it\}\), take the straight segment from \(z_b\) to \(p-h\), the counterclockwise circle of radius \(h\) about \(p\), and the return segment. Denote these three based loops by \(\gamma_0,\gamma_+,\gamma_-\). The product of loops means traversal from left to right. Write \(G_\gamma\) for the deck transformation that takes the initial local lift \(\tau\) to the branch reached by analytic continuation around \(\gamma\). This transformation is the monodromy of the lift around \(\gamma\). Continuation along the product \(\gamma\nu\) gives \(G_\gamma G_\nu\tau\). Since all three prescribed points are simple, each image loop winds once positively about the corresponding puncture. Use a deck transformation to normalize the first monodromy to \(A\). The second is conjugate to \(B\), where
\[
 A=\begin{pmatrix}1&2\\0&1\end{pmatrix},\qquad
 B=\begin{pmatrix}1&0\\-2&1\end{pmatrix}.
\]
These matrices represent the positive cusp translations for \(0\) and \(1\), respectively, in the coordinates \(\tau\) and \(-1/\tau\).

For a holomorphic map from an annulus \(1<|\xi|<R_{\mathrm a}\) to \(X\), take a counterclockwise circle in the annulus as a positive core loop. Here \(R_{\mathrm a}\) is the outer radius. Let \(T\) be the deck transformation of this loop, and let \(\tr T\) be the trace, the sum of the diagonal entries of its determinant-one matrix. The absolute trace is independent of the sign of the matrix. If \(|\tr T|>2\), \(T\) is hyperbolic. The Schwarz--Pick inequality for the hyperbolic metrics of curvature \(-1\) gives \cite[Section~2]{BE}
\begin{equation}\label{eq:length}
 2\arcosh(|\tr T|/2)\le\frac{2\pi^2}{\log R_{\mathrm a}}.
\end{equation}
Here \(\arcosh\) is the inverse hyperbolic cosine, and the left side is the hyperbolic translation length of \(T\). The right side is the shortest hyperbolic length of a core loop. In logarithmic coordinates, the annulus unfolds to a strip of width \(\log R_{\mathrm a}\), on whose center line the density is \(\pi/\log R_{\mathrm a}\) and the core translates by \(2\pi i\).
For \(R_{\mathrm a}=17,19\) we have \(\pi^2<\log R_{\mathrm a}\,\arcosh R_{\mathrm a}\). Indeed, \(\log17>2.83,\ \arcosh17>3.52,\ \log19>2.9,\ \arcosh19>3.6\) follow from \(\arcosh R_{\mathrm a}=\log(R_{\mathrm a}+\sqrt{R_{\mathrm a}^2-1})\) and positive-term truncation of \(\log x=2\sum_{j\ge0}v^{2j+1}/(2j+1)\), where \(v=(x-1)/(x+1)\). Combining these lower bounds with \(\pi<22/7\) gives the claimed inequalities. Thus \eqref{eq:length} excludes \(|\tr T|\ge2R_{\mathrm a}\), yielding strict trace bounds of \(34\) and \(38\), respectively.

\medskip
\noindent\textit{The first annulus determines the transformation around \(it\).}
Put \(\zeta=iz/t\), so that \(0,it,-it\) map to \(0,-1,1\). Cut the disk \(|\zeta|<10\) along \([1,10)\). The map
\begin{equation}\label{eq:annulus-map}
 u=\frac{9(10+\zeta)}{11(10-\zeta)},\qquad
 v=\sqrt{\frac{u^2}{1-u^2}}
\end{equation}
takes the cut disk conformally onto the right half-plane, with the square root positive for \(0<u<1\). The first step gives the right half-plane slit along \([1,\infty)\). Squaring gives the plane slit along \((-\infty,0]\) and \([1,\infty)\), and the fractional linear map and square root give the asserted half-plane.
The segment \([-1,0]\) maps to a positive real segment \([v_-,v_+]\). Evaluating \eqref{eq:annulus-map} at the endpoints gives \(\ell=\sqrt{v_+/v_-}=(202/81)^{1/4}<81/64\). The change of coordinate
\[
 \chi=\frac{v-\sqrt{v_-v_+}}{v+\sqrt{v_-v_+}}
\]
maps the half-plane to the unit disk and the segment to \([-\kappa,\kappa]\), with \(\kappa=(\ell-1)/(\ell+1)<17/145\).
Thus \(\chi=\kappa(\xi+\xi^{-1})/2\) maps \(1<|\xi|<17\) injectively into the unit disk and avoids the segment, since \(\kappa(17+17^{-1})/2<1\). Inverting the preceding maps embeds the annulus in the three-puncture source disk.

Its core comes from a positively oriented ellipse around \([-\kappa,\kappa]\). In the cut \(\zeta\)-disk it encloses just \(0\) and \(-1\). The boundary of a thin neighborhood of the two small circles and their tails identifies the core with \(\gamma_0\gamma_+\) up to conjugacy and cyclic permutation. If the transformation of \(\gamma_+\) is \(MBM^{-1}\) with \(M\in\Gamma(2)\) and entries \(m_{jk}\) in row \(j\) and column \(k\), matrix multiplication gives \(\tr(AMBM^{-1})=2-4m_{22}^2\). The entry \(m_{22}\) is odd, so the strict trace bound \(34\) forces \(m_{22}=\pm1\). The fixed point \(m_{12}/m_{22}\) of \(MBM^{-1}\) is then an even integer. Simultaneous conjugation by the translation \(A^{-m_{12}/(2m_{22})}\) preserves \(A\) and moves that point to \(0\). The new conjugating matrix is lower triangular, with both diagonal entries \(1\) or both \(-1\), and commutes with \(B\). The transformations around \(0,it\) are therefore \(A,B\).

\medskip
\noindent\textit{The second annulus determines the fixed line of the reflection.}

The real symmetry of \(f\) lifts at the base point to an antiholomorphic reflection \(\mathcal R=\gamma\circ J_0\) with \(J_0(\tau)=-\bar\tau\) and \(\gamma\in\Gamma(2)/\{\pm I\}\). Choose it to fix \(\tau(z_b)\). Reflection sends \(\gamma_0\) to its inverse. Since \(J_0 A J_0=A^{-1}\), the transformation \(\gamma\) commutes with \(A\) and therefore fixes \(\infty\). The stabilizer of \(\infty\) in the deck group consists of the \(A^n\) with \(n\in\mathbb Z\). Thus
\[
 \mathcal R(\tau)=2n-\bar\tau,\qquad \re\tau(z_b)=n.
\]
As \(f(z_b)<0\), the identities \eqref{eq:lambda} force \(n\) to be odd. Reflection exchanges \(\gamma_+\) with \(\gamma_-^{-1}\). Using \(J_0 B J_0=B^{-1}\), the monodromy of \(\gamma_-\) is \(A^n B A^{-n}\).

The counterclockwise boundary of a thin regular neighborhood of the three tails and circles, traced from the left of the base point, visits the lower, middle, and upper branches in that order. Thus the positive outer loop has word \(\gamma_-\gamma_0\gamma_+\) up to conjugacy. The annulus
\[
 z=\tfrac{it}{2}(\xi+\xi^{-1}),\qquad 1<|\xi|<19,
\]
is univalent, avoids the three punctures, and lies in the unit disk because \((19+19^{-1})/20=181/190\). Its core is a positive ellipse around the three punctures, so it represents that outer loop. Multiplication and invariance of trace under cyclic permutation give
\[
 \tr(ABA^n B A^{-n})=-(16n(n-1)+6).
\]
For odd \(n\ne1\) the absolute trace is at least \(38\). The strict trace bound \(38\) forces \(n=1\), so the reflection fixes the line \(\re\tau=1\).

\medskip
\noindent\textit{The common branch.}
Continue the normalized lift from the base point just into \(\Om\), using the upper tail to identify the positive loop around \(it\). The positive generator of the fundamental group of \(\Om\setminus\{it\}\) has monodromy \(B\), which translates \(-1/\tau\) by \(2\). Thus \(\phi=-\exp(-\pi i/\tau)\) is single-valued there. Boundedness gives a holomorphic extension across \(it\), and the maximum modulus principle keeps its range in \(\D\). By continuity, \(U(\phi(it))=f(it)=1\), so \(\phi(it)=0\).

The reflection with \(n=1\) fixes the branch along the negative real segment on the line \(\re\tau=1\). Its quotient image therefore lies on \(C\). Near the source point \(0\), the local cusp coordinate \(q=\Lambda^{-1}(f(z))\), with the inverse branch taking \(0\) to \(0\), is real-symmetric with \(q'(0)>0\). In a sufficiently small upper half-neighborhood it has positive imaginary part, by the Taylor estimate used in Proposition~\ref{pr:candidate}. Along an upper semicircle from the negative to the positive segment, the endpoint arguments of this \(q\)-coordinate are \(\pi\) and \(0\). The \(\tau\)-branch therefore has real part \(0\) on the adjacent positive segment. Since \(0<f(x)<1\) on \((0,1)\), continuation along this segment stays on the line \(\re\tau=0\). Its quotient image lies in \((-1,0)\). Local inverses of the covering give analytic continuation across both open segments.
\end{proof}

\subsection{Exclusion and the equality case}\label{sec:exclude}

Use the map
\[
 K(w)=\log\!\left(\frac{1+w}{1-w}\right)
\]
with the branch on the right half-plane that vanishes at \(1\). It is injective on the disk, and \(|\im K|<\pi/2\) there. For \(w=e^{-g_C(\theta)+i\theta}\in C\) with \(g_C(\theta)=\sqrt{\theta(\pi-\theta)}\), direct evaluation gives
\begin{equation}
 \re K(w)=\atanh\!\left(\frac{\cos\theta}{\cosh g_C(\theta)}\right),\qquad
 \im K(w)=\arctan\!\left(\frac{\sin\theta}{\sinh g_C(\theta)}\right).
 \label{eq:graph}
\end{equation}
The real part decreases strictly from \(+\infty\) to \(-\infty\), since \(g_C'(\theta)\cos\theta\ge0\) and \(\sin\theta>0\). The imaginary part is positive and tends to zero at both ends. On \(K(C)\), each real part therefore determines a unique imaginary part, as required by Lemma~\ref{lem:primitive}. The imaginary part is zero on \((-1,0)\) and has the interior limit \(-\pi/2\) on \(L\).

Suppose \(0<t\le\min(t_*,1/10)\) is feasible for the disk problem and normalize \(f'(0)>0\). Use Proposition~\ref{pr:lift} to obtain \(\phi\), and put \(\rho=t/t_*\). Choose a conformal map \(\eta:\D\to\Om\) with \(\eta(0)=it_*\). The map extends homeomorphically to the boundary. On \(\D\) define
\begin{equation}
 W_0=K\circ\Phi_*\circ\eta,\qquad
 W_1(\zeta)=K\bigl(\phi(\rho\,\eta(\zeta))\bigr).
 \label{eq:compare-lifts}
\end{equation}
Both vanish at \(0\). On the arc corresponding to the negative real segment, the two boundary values lie on the graph \eqref{eq:graph}. On the arc corresponding to the positive segment both imaginary parts are zero. On the remaining arc the candidate has imaginary part \(-\pi/2\), so the competitor's radial boundary values are at least as large almost everywhere. If \(\rho<1\), the scaled open semicircle lies inside \(\Om\), and the competitor's imaginary part is strictly greater at every point of this arc.

At the two ends of the negative segment the candidate's imaginary part tends to zero, and the competitor's imaginary part on this segment is nonnegative. Along the positive segment the difference vanishes, and along the outer arc it is nonnegative almost everywhere. These facts give the two essential lower limits at each of the three junctions, as required by Lemma~\ref{lem:primitive}. The slit-tip image on the open positive segment also has imaginary part zero. We conclude that \(W_1=W_0\). The strict inequality on the outer arc excludes \(\rho<1\). For \(\rho=1\), the injectivity of \(K\) gives \(\phi=\Phi_*\), and the identity theorem gives \(f=f_*\) on \(\D\).

The certified rational controller at \(\delta=98236947/10^8\) gives, by \eqref{eq:parameter}, a feasible disk parameter \(t_w<1/10\). The certificate is available in the repository linked in Section~\ref{sec:conclusion}. If \(t_w<t_*\), the preceding comparison would exclude that interpolant. Thus \(t_*\le t_w<1/10\), so the comparison covers every \(0<t\le t_*\). Existence at \(t_*\) follows from Proposition~\ref{pr:candidate}. For \(t_*<t<1\), the function \(z\mapsto f_*((t_*/t)z)\) is a disk interpolant at \(t\). This proves Theorem~\ref{th:disk}.

\section{Finite-order controllers}\label{sec:finite}

By \eqref{eq:cayley}, feasibility at \(\delta\) implies disk feasibility at the parameter in \eqref{eq:parameter}. Theorem~\ref{th:disk} then gives \(t\ge t_*\), hence \(\delta\le\delta_*\) by \eqref{eq:mass}.

\subsection{The endpoint is unattained}

Suppose there is a feasible controller at \(\delta_*\). Normalize the corresponding rational function \(R\) from \eqref{eq:cayley}, precomposing with \(z\mapsto-z\) if needed. By Theorem~\ref{th:disk} the normalized rational function agrees with \(f_*\) on the disk. At the boundary point \(z=i\) this rational function is holomorphic and takes a value in \(X\). The Cayley coordinates of \(i,-i\) are \(i,-i\), respectively. The polynomials \(x,y,Z\) have no imaginary-axis zeros, and \(a_{\delta_*}(\pm i)\) and \(b(\pm i)\) are both nonzero.

Take a small disk in \(X\) about this boundary value that is evenly covered by \(U:\D\setminus\{0\}\to X\). This means its preimage is a union of disjoint open sets, called sheets, each mapped biholomorphically onto the disk by \(U\). A small neighborhood of \(i\) maps into that disk, and its intersection with \(\Om\) is connected. Since \(f_*=U\circ\Phi_*\) in \(\Om\), the lift \(\Phi_*\) on this intersection lies in a single sheet. The corresponding local inverse extends \(\Phi_*\) to \(i\) with value in \(\D\setminus\{0\}\). But \(\Phi_*\) maps that open boundary arc into \(L\), giving a limit of modulus one. This contradiction excludes the endpoint.

\subsection{Construction below the threshold}

Fix \(0<\delta<\delta_*\) and put \(\sigma=\sqrt{(1-\delta)/(1+\delta)}\) and \(\rho=t_*/\sigma<1\). Then \(F(z)=f_*(\rho z)\) is holomorphic on \(|z|<1/\rho\), with just one simple zero at \(0\) and two simple \(1\)-points at \(\pm i\sigma\). Define
\[
 P_\sigma(z)=1+z^2/\sigma^2,\qquad
 D_\infty(z)=\frac{P_\sigma(z)}{1-F(z)},\qquad
 A_\infty(z)=\frac{D_\infty(z)-P_\sigma(z)}{z}=\frac{D_\infty(z)F(z)}{z}.
\]
The numerator \(P_\sigma\) and denominator \(1-F\) have simple zeros at \(\pm i\sigma\), so their quotient extends there with nonzero values. Likewise \(F(z)/z\) extends at \(0\) with nonzero value \(F'(0)\). Thus \(D_\infty,A_\infty\) are holomorphic, zero-free and real-symmetric on the larger disk, with \(D_\infty(0)=1\). Let \(D_m\) be the Taylor polynomial for \(D_\infty\) through degree \(m\ge2\) and set \(A_m=(D_m-P_\sigma)/z\). This is a real polynomial because \(D_m(0)=P_\sigma(0)=1\). On the closed unit disk, \(D_m\) tends uniformly to \(D_\infty\). The function \((D_m-D_\infty)/z\) is holomorphic across \(0\), and on the unit circle its modulus equals \(|D_m-D_\infty|\). The maximum modulus principle therefore gives uniform convergence of \(A_m\) to \(A_\infty\) on the closed disk. Since both limits have positive minimum modulus there, \(D_m\) and \(A_m\) have no zeros in the closed disk for all sufficiently large \(m\).

Choose such an \(m\), take \(n\ge\max(1,\deg A_m)\), and use the Cayley transform \(z=(s-1)/(s+1)\) to define
\begin{equation}
 x(s)=\frac{2(s+1)^n}{1-\delta},\qquad y(s)=(s+1)^n A_m(z).
 \label{eq:controller}
\end{equation}
Using \((s+1)^2P_\sigma(z)=2a_\delta(s)/(1-\delta)\) and \((s+1)^2z=b(s)\), we obtain
\[
 Z(s)=(s+1)^{n+2}\bigl(P_\sigma(z)+zA_m(z)\bigr)
     =(s+1)^{n+2}D_m(z).
\]
The bounds \(\deg A_m\le n\) and \(\deg D_m\le n+1\) make these real polynomials. Since \(A_m(1)\ne0\), we have \(\deg y=\deg x=n\). For \(\re s\ge0\), the Cayley image has \(|z|\le1\) and \(s+1\ne0\). Hence \(x,y,Z\) have no zeros in the closed right half-plane. This proves feasibility for every \(0<\delta<\delta_*\). Together with the endpoint argument, the parameter-domain check in Section~\ref{sec:intro}, and Proposition~\ref{pr:balance}, this completes the proof of Theorem~\ref{th:main}.

\section{Two results for fixed-degree controllers}\label{sec:fixed-degree}

Let \(H^i\) denote the set of real strictly Hurwitz polynomials of degree \(i\), and let \(MH^n\) denote its monic counterpart at degree \(n\). Thus \(H^0\) consists of the nonzero real constants. For \(0\le i\le n\), let \(S_{n,i}\) be the set of positive \(\delta\) for which \eqref{eq:poly} holds with \(x\in MH^n\), \(y\in H^i\), and \(Z\) strictly Hurwitz. Write \(\delta_n^*=\sup S_{n,0}\). Theorem~\ref{th:main} places all these sets inside \((0,\delta_*)\).

The set \(S_{n,0}\) is nonempty for every \(n\ge1\). Indeed, take \(x=(s+r)^n\), with \(r>0\) large enough that \(\im x(i)>0\). At \(\delta=0\), the simple root \(i\) of \(a_0x\) moves with derivative \(-i/x(i)\) as the constant \(y\) increases from zero. This derivative has negative real part. The conjugate root moves to the left as well, and the roots near \(-r\) stay there by continuity. A sufficiently small positive \(y\) makes \(Z\) strictly Hurwitz at \(\delta=0\), and stability persists for small positive \(\delta\).

For \(\delta\in S_{n,0}\), put \(\theta=\arccos\delta\) and \(\zeta=e^{i\theta}\). Then \(a_\delta(\zeta)=0\) and \(\zeta^{-1}Z(\zeta)=2iy\sin\theta\). Writing \(Z(s)=\sum_{j=0}^{n+2}z_js^j\), we obtain
\begin{equation}\label{eq:fixed-balance}
 \sum_{j=0}^{n+2}z_j\cos((j-1)\theta)=0.
\end{equation}
Here \(z_{n+2}=1\) and every \(z_j>0\), since a real monic strictly Hurwitz polynomial factors into real linear and quadratic factors with positive coefficients.

\subsection{The effect of numerator degree at order six}

For a fixed denominator degree \(n\), consider whether an upper bound \(\delta_b^*\) on the feasible parameters with \(x\in MH^n,\ y\in H^0\) also holds with \(y\in H^i\) for \(1\le i\le n\) \cite[Appendix~C]{HeWY}. The following example gives a counterexample at \(n=i=6\).

\begin{example}[A sixth-order counterexample]\label{ex:fixed-six}
Put \(\delta_b^*=\sqrt{11/12}\). Every \(\delta\in S_{6,0}\) satisfies \(\delta<\delta_b^*\), and some \(\delta\in S_{6,6}\) exceeds \(0.96295>\delta_b^*\). In particular,
\[
 \delta_6^*\le\delta_b^*<0.96295<\sup S_{6,6}.
\]
\end{example}

\begin{proof}[Proof outline]
First exclude a constant numerator. If \(\delta\ge\sqrt{11/12}\), then \(0<\theta<\pi/10\). For a strictly Hurwitz polynomial \(x\) of degree six, the continuous argument of \(x(i\omega)\) increases from \(0\) to \(3\pi\) as \(\omega\) increases from zero to infinity, because the argument of each root factor increases strictly. At argument \(3\pi/2\), the imaginary part is negative, so \(x_1-\omega^2 x_3+\omega^4x_5<0\), where \(x_k\) is the coefficient of \(s^k\). Coefficient comparison in \eqref{eq:poly} expresses \(x_1,x_3,x_5\) as linear functions of the closed-loop coefficients \(z_0,\ldots,z_8\). Write \(C_j=\cos((j-1)\theta)\). The nonnegative coefficient vectors satisfying \eqref{eq:fixed-balance} are generated by coordinate vectors for which \(C_j=0\) and by vectors \(-C_j e_i+C_i e_j\) for \(C_i>0>C_j\), where \(e_i,e_j\) are coordinate vectors. Only \(j=7,8\) can have nonpositive weights. Substituting each generator into the coefficient recurrence gives odd coefficients \((v_1,v_3,v_5)\) with \(v_5>0\) and \(4v_1v_5-v_3^2\ge0\). Thus every generator makes \(v_1-\omega^2 v_3+\omega^4v_5\) nonnegative. Nonnegative combinations preserve this inequality, contradicting the strict Hurwitz sign condition.

For the sixth-degree numerator, let \(A\) be the root of \(1+2A-2A^2-2A^3=0\) in \((0.854637,0.854638)\), and set \(d=\sqrt{(1+A)/2}\) and \(k=3A^2-1\). Expansion gives
\[
 \begin{gathered}
 X_0(s)=(s^2+2ds+1)(s^2+A)^2,\qquad
 Y_0(s)=k(s^2+A^2/k),\\
 a_dX_0+bY_0=s^8.
 \end{gathered}
\]
Set \(h=10^{-4}\) and \(K=1+h^2-2dh\), and apply the fractional linear substitution \(s\mapsto(s+h)/(1+hs)\) as in Charles and Boston \cite{CB}. It maps the closed right half-plane into the open right half-plane. Clearing denominators and multiplying the transformed \(X_0,Y_0\) by \(K\) and \(1-h^2\), respectively, gives two sixth-degree strictly Hurwitz polynomials with closed loop \((s+h)^8\). A common rescaling makes the denominator monic. The transformed parameter is \(\delta'=[d(1+h^2)-2h]/K\). Since \(d>0.962973\) and this expression increases with \(d\), we obtain \(\delta'>0.96295>\sqrt{11/12}\).
\end{proof}

\subsection{The continued-fraction bound for constant numerators}

For constant numerators, consider whether the largest positive real zero \(\delta_n\) of the following continued fraction bounds the supremum \(\delta_n^*\) \cite[Section~2.3]{HeWY}. Define
\[
 f_0(\delta)=-\delta,\qquad f_n(\delta)=-\frac{1}{f_{n-1}(\delta)}-2\delta\quad(n\ge1),
\]
as rational functions of \(\delta\), and let \(\delta_n\) be the largest positive real zero of \(f_n\).

\begin{theorem}[Sharpness of the continued-fraction bound]\label{th:fixed-bound}
For every \(n\ge1\),
\[
 \delta_n=\sin\!\left(\frac{(n+2)\pi}{2(n+1)}\right)
         =\cos\!\left(\frac{\pi}{2(n+1)}\right),\qquad
 \delta_n^*\le\delta_n.
\]
Every \(\delta\in S_{n,0}\) is strictly below \(\delta_n\). The supremum satisfies
\[
 \delta_n^*=\delta_n\quad(1\le n\le4),\qquad
 \delta_n^*<\delta_n\quad(n\ge5).
\]
\end{theorem}

\begin{proof}[Proof outline]
Let \(T_j\) be the Chebyshev polynomials of the first kind, with \(T_j(\cos\theta)=\cos(j\theta)\). The continued fraction reduces to \(f_n(\delta)=-T_{n+1}(\delta)/T_n(\delta)\). Consecutive polynomials have no common root, so the largest positive zero is \(\cos(\pi/[2(n+1)])\). If \(\delta\ge\delta_n\) were feasible, then \(\theta\le\pi/[2(n+1)]\) would make every term on the left of \eqref{eq:fixed-balance} nonnegative, with the constant term strictly positive.

Set \(m=n+2\), \(\theta_0=\pi/[2(n+1)]\), and \(c_0=1/(2\sin\theta_0)\). If feasible parameters approach \(\delta_n\), every weight in \eqref{eq:fixed-balance} except the leading one tends to a positive number, while the leading weight tends to zero. Hence \(Z\) tends coefficientwise to \(s^m\), and taking the imaginary part of \(\zeta^{-1}Z(\zeta)=2iy\sin\theta\) shows that \(y\) tends to \(c_0\). The denominator \(x\) must therefore tend to
\[
 \frac{s^m-c_0(s^2-1)}{a_{\delta_n}(s)}.
\]
For \(n\ge5\), the numerator of the limit has a root with argument in \((2\pi/(m-2),3\pi/m)\). Substituting \(s=re^{i\varphi}\) into the numerator and eliminating \(r\) gives the continuous function
\[
 \sin^m(m\varphi)-c_0^2\sin^2(2\varphi)\sin^{m-2}((m-2)\varphi).
\]
It has opposite signs at the two ends of the interval. At an interior zero, take the positive radius \(r\) with \(r^2=\frac{\sin(m\varphi)}{\sin((m-2)\varphi)}\). The resulting root lies in the open right half-plane and differs from both roots of \(a_{\delta_n}\). Continuity of roots excludes strictly Hurwitz denominators approaching this limit.

For \(1\le n\le4\), take \(\theta\) slightly above \(\theta_0\), put \(\varphi=(\theta+\pi/2)/m\), and choose
\[
 \varepsilon=\frac{\sin(\theta-\varphi)}{\sin\varphi},\qquad
 y=\frac{(\sin\theta/\sin\varphi)^m}{2\sin\theta},\qquad Z=(s+\varepsilon)^m.
\]
Then \(Z-yb\) vanishes at \(e^{\pm i\theta}\), so division by \(a_{\cos\theta}\) gives a monic polynomial \(x\). In degrees one through three the limiting polynomial is strictly Hurwitz. In degree four the limit has a pair of imaginary-axis roots. For the perturbed polynomial \(x=s^4+x_3s^3+x_2s^2+x_1s+x_0\), coefficient comparison shows that the Hurwitz determinant \(x_3x_2x_1-x_1^2-x_3^2x_0\), divided by \(\varepsilon\), has a positive limit. Its coefficients stay positive. The construction thus gives strictly feasible parameters \(\cos\theta\) arbitrarily close to \(\delta_n\) in all four degrees.
\end{proof}

\section{Conclusion}\label{sec:conclusion}

Using a balanced slit domain and a comparison of harmonic primitives, we prove the exact stabilizability threshold formula for the Belgian Chocolate Problem and classify all disk interpolants at the critical parameter.

The key step is to turn the boundary markings of a competing lift into a comparison of harmonic primitives. Two explicit annuli determine the monodromy and hence the boundary markings. On the graph arc, equality of the tangential derivatives at a boundary minimum forces equality of the radial derivatives. The one-sided bounds on the remaining arcs and the Poisson-kernel estimate then force the normalized lifts to agree.

Let \(d_*\) be the harmonic measure at the origin of the lower semicircle of the balanced slit domain. The feasible disk parameters form \([t_*,1)\), where \(t_*=\tan(\pi d_*/4)\). At \(t_*\), there are exactly two interpolants. One is \(f_*\), with \(f_*'(0)>0\). The other is obtained by replacing the argument of \(f_*\) by its negative. The feasible finite-order controller parameters form \((0,\delta_*)\), where \(\delta_*=\cos(\pi d_*/2)\).

The extremal interpolant explains how these endpoints arise. Its lift approaches the unit circle at \(z=i\). A finite-order controller at \(\delta_*\) would force the same lift to extend there with a value inside the disk, which excludes attainment. For each \(0<\delta<\delta_*\), rescaling the extremal allows uniform polynomial approximation to preserve the required nonvanishing conditions and yield a feasible controller.

The fixed-degree results show how restrictions on the controller affect stabilizability. At denominator degree six, allowing a nonconstant numerator strictly raises the supremum of feasible parameters. For constant numerators, the continued-fraction bound \(\delta_n\) is sharp at degrees one through four. From degree five onward, \(\delta_n^*<\delta_n\).

We verified the threshold theorems end to end in Lean~4 \cite{Lean4}, using Mathlib \cite{Mathlib} and Tau Ceti \cite{TauCeti}. Starting from the concrete slit domains and the two feasibility definitions, the formalization establishes the exact threshold formulas, both feasible intervals, uniqueness of the normalized critical interpolant, and nonattainment by finite-order controllers. Lean~4 also verifies the six-degree separation and the full sharpness classification. The global threshold development and its reproduction instructions are available at \href{https://github.com/Jxxxx-2001/belgian-chocolate-lean}{\textcolor{blue}{\nolinkurl{https://github.com/Jxxxx-2001/belgian-chocolate-lean}}}.

\section*{AI Declaration}

The authors led the formulation of the research question, set the research objectives, designed the overall proof framework, and organized and interpreted the results. GPT 6.0 Astra obtained the main result on the exact stabilizability threshold in October 2026. Coding agents using GPT-6.2 Astra and GPT-6.2 Sol subsequently assisted with the Lean~4 formalization. The authors wrote the manuscript and used AI to improve its grammar, spelling, and readability. The authors take full responsibility for the content of the article.

\renewcommand{\bibfont}{\small}
\setlength{\bibsep}{3pt plus .5pt minus .5pt}
\bibliographystyle{unsrtnat}
\bibliography{references}

@article{Youla,
  author = {Youla, D. C. and Bongiorno, Jr., J. J. and Lu, C. N.},
  title = {Single-loop feedback-stabilization of linear multivariable dynamical plants},
  journal = {Automatica},
  year = {1974},
  volume = {10},
  number = {2},
  pages = {159--173},
  doi = {10.1016/0005-1098(74)90021-1}
}

@book{Blondel,
  author = {Blondel, Vincent D.},
  title = {Simultaneous Stabilization of Linear Systems},
  series = {Lecture Notes in Control and Information Sciences},
  volume = {191},
  publisher = {Springer},
  address = {Berlin},
  year = {1994},
  doi = {10.1007/3-540-19862-8}
}

@article{GuanHeWangYu,
  author = {Guan, Qiang and He, Guan-nan and Wang, Long and Yu, Wen-sheng},
  title = {Simultaneous stabilization of linear systems},
  journal = {Control Theory \& Applications},
  year = {2011},
  volume = {28},
  number = {1},
  pages = {1--12},
  note = {In Chinese}
}

@article{LiWangYu,
  author = {Li, Wang and Wang, Long and Yu, Wensheng},
  title = {Some open problems on simultaneous stabilization of linear systems},
  journal = {Journal of Systems Science and Complexity},
  year = {2016},
  volume = {29},
  number = {2},
  pages = {289--299},
  doi = {10.1007/s11424-015-4182-1}
}

@article{BE,
  author = {Bergweiler, Walter and Eremenko, Alexandre},
  title = {Gol'dberg's constants},
  journal = {Journal d'Analyse Math\'ematique},
  year = {2013},
  volume = {119},
  pages = {365--402},
  doi = {10.1007/s11854-013-0012-3}
}

@article{Patel,
  author = {Patel, Vijay V. and Deodhare, Girish and Viswanath, T.},
  title = {Some applications of randomized algorithms for control system design},
  journal = {Automatica},
  year = {2002},
  volume = {38},
  number = {12},
  pages = {2085--2092},
  doi = {10.1016/S0005-1098(02)00112-7}
}

@article{BHLO,
  author = {Burke, James V. and Henrion, Didier and Lewis, Adrian S. and Overton, Michael L.},
  title = {Stabilization via nonsmooth, nonconvex optimization},
  journal = {IEEE Transactions on Automatic Control},
  year = {2006},
  volume = {51},
  number = {11},
  pages = {1760--1769},
  doi = {10.1109/TAC.2006.884944}
}

@article{ChangSahinidis,
  author = {Chang, YoungJung and Sahinidis, Nikolaos V.},
  title = {Global optimization in stabilizing controller design},
  journal = {Journal of Global Optimization},
  year = {2007},
  volume = {38},
  number = {4},
  pages = {509--526},
  doi = {10.1007/s10898-006-9092-2}
}

@article{Batra2007,
  author = {Batra, Prashant},
  title = {On {Goldberg}'s Constant {$A_2$}},
  journal = {Computational Methods and Function Theory},
  year = {2007},
  volume = {7},
  number = {1},
  pages = {33--41},
  doi = {10.1007/BF03321629}
}

@article{HeWY,
  author = {He, Guannan and Wang, Long and Yu, Wensheng},
  title = {French champagne and {Belgian} chocolate problems in simultaneous stabilization of linear systems},
  journal = {IFAC Proceedings Volumes},
  year = {2008},
  volume = {41},
  number = {2},
  pages = {1105--1110},
  doi = {10.3182/20080706-5-KR-1001.00191}
}

@inproceedings{LiYu,
  author = {Li, Wang and Yu, Wensheng},
  title = {A sharper parameter bound for {Belgian} chocolate stabilization problem},
  booktitle = {Proceedings of the 26th Chinese Control and Decision Conference},
  year = {2014},
  pages = {459--462},
  doi = {10.1109/CCDC.2014.6852192}
}

@article{CB,
  author = {Charles, Zachary and Boston, Nigel},
  title = {Exploiting algebraic structure in global optimization and the {Belgian} chocolate problem},
  journal = {Journal of Global Optimization},
  year = {2018},
  volume = {72},
  number = {2},
  pages = {241--254},
  doi = {10.1007/s10898-018-0659-5}
}

@misc{ChengBCP,
  author = {Cheng, Xuandi},
  title = {Certified lower bound for the {Belgian Chocolate Problem}},
  year = {2026},
  howpublished = {GitHub technical note and verification repository},
  note = {Revision \texttt{362cf25a}},
  url = {https://github.com/chengxuandi/belgian-chocolate-certified-lower-bound/tree/362cf25a04ebd0468b05153008ba72dca7f09ef6}
}

@misc{Eremenko,
  author = {Eremenko, Alexandre},
  title = {Simultaneous stabilization, avoidance and {Goldberg}'s constants},
  year = {2012},
  note = {arXiv:1208.0778},
  url = {https://arxiv.org/abs/1208.0778}
}

@misc{Liu,
  author = {Liu, Muhan},
  title = {An analytic solution to the generalized {Belgian} chocolate problem},
  year = {2026},
  note = {Preprint, version 1},
  url = {https://zenodo.org/records/24108157}
}

@book{Pommerenke,
  author = {Pommerenke, Christian},
  title = {Boundary Behaviour of Conformal Maps},
  series = {Grundlehren der mathematischen Wissenschaften},
  volume = {299},
  publisher = {Springer},
  address = {Berlin},
  year = {1992},
  doi = {10.1007/978-3-662-02770-7}
}

@article{BRY,
  author = {Binder, Ilia and Rojas, Cristobal and Yampolsky, Michael},
  title = {Carath\'eodory convergence and harmonic measure},
  journal = {Potential Analysis},
  year = {2019},
  volume = {51},
  number = {4},
  pages = {499--509},
  doi = {10.1007/s11118-018-9721-7}
}

@misc{DLMF,
  author = {{National Institute of Standards and Technology}},
  title = {{NIST Digital Library of Mathematical Functions}},
  year = {2026},
  note = {Version 1.2.8, Sections 23.15, 23.17, and 23.18},
  url = {https://dlmf.nist.gov/}
}

@book{HFT,
  author = {Axler, Sheldon and Bourdon, Paul and Ramey, Wade},
  title = {Harmonic Function Theory},
  edition = {2nd},
  series = {Graduate Texts in Mathematics},
  volume = {137},
  publisher = {Springer},
  address = {New York},
  year = {2001},
  url = {https://www.axler.net/HFT.pdf}
}

@inproceedings{Lean4,
  author = {de Moura, Leonardo and Ullrich, Sebastian},
  title = {The {Lean 4} theorem prover and programming language},
  booktitle = {Automated Deduction -- CADE 28},
  series = {Lecture Notes in Computer Science},
  volume = {12699},
  publisher = {Springer},
  year = {2021},
  pages = {625--635},
  doi = {10.1007/978-3-030-79876-5_37}
}

@inproceedings{Mathlib,
  author = {{The mathlib Community}},
  title = {The {Lean} mathematical library},
  booktitle = {Proceedings of the 9th ACM SIGPLAN International Conference on Certified Programs and Proofs},
  publisher = {ACM},
  year = {2020},
  pages = {367--381},
  doi = {10.1145/3372885.3373824}
}

@misc{TauCeti,
  author = {{Tau Ceti Project}},
  title = {{Tau Ceti}: a {Lean} library downstream of {Mathlib}},
  year = {2026},
  note = {Software, revision \texttt{ed6c0aa3}},
  url = {https://github.com/TauCetiProject/TauCeti/tree/ed6c0aa3acc0a3cb56901dc406167cdbad593f59}
}

@misc{LiuLean,
  author = {Liu, Muhan},
  title = {{Lean} formalization of the {Belgian} chocolate problem},
  year = {2026},
  howpublished = {Source code accompanying the preprint},
  url = {https://github.com/Panda-China666/belgian-chocolate-lean/blob/39ea4c0e2930ded15e3e6fcc0bbcbebbf6d6b89e/BelgianChocolate.lean#L3913}
}

\end{document}